\documentclass[11pt]{elsarticle}
\usepackage[a4paper,margin=3cm]{geometry}
\usepackage{amsmath,amssymb}
\usepackage{amsthm}
\usepackage{mathtools}
\usepackage{breqn}
\usepackage{hyperref}
\usepackage{cleveref}

\theoremstyle{definition}
\newtheorem{definition}{Definition}[section]

\theoremstyle{plain}

\theoremstyle{remark}
\newtheorem{remark}[definition]{Remark}

\theoremstyle{definition}

\makeatletter
\newcommand*{\defeq}{\mathrel{\rlap{%
					 \raisebox{0.3ex}{$\m@th\cdot$}}%
					 \raisebox{-0.3ex}{$\m@th\cdot$}}%
					 =}
\makeatother

\makeatletter
\newcommand*{\eqdef}{=\mathrel{\rlap{%
					 \raisebox{0.3ex}{$\m@th\cdot$}}%
					 \raisebox{-0.3ex}{$\m@th\cdot$}}%
					 }
\makeatother

\makeatletter
\newcommand{\proofstep}[1]{%
	\par
	\addvspace{\medskipamount}%
	\noindent\textit{#1\@addpunct{.}}\enspace\ignorespaces
}
\makeatother

\DeclarePairedDelimiter\abs{\lvert}{\rvert}%
\DeclarePairedDelimiter\norm{\lVert}{\rVert}%
\DeclarePairedDelimiter\croc{\langle}{\rangle}%

\makeatletter
\let\oldabs\abs
\def\abs{\@ifstar{\oldabs}{\oldabs*}}
\let\oldnorm\norm
\def\norm{\@ifstar{\oldnorm}{\oldnorm*}}
\let\oldcroc\croc
\def\croc{\@ifstar{\oldcroc}{\oldcroc*}}
\makeatother

\newcommand{\ee}{\operatorname{e}}

\begin{document}

\begin{frontmatter}
	
\title{Analytic umbral transmutations and Bessel moments}

\author{Roberto Ricci\corref{cor1}}
\ead{roberto.ricci@enea.it}

\author{Giuseppe Dattoli}

\cortext[cor1]{Corresponding author}

\address{ENEA, Nuclear Department, Frascati Research Center, Via E. Fermi 45, 00044 Frascati (Rome), Italy}

\begin{abstract}
	We develop an analytic umbral approach to Bessel moments, using them as a concrete testbed justifying the passage from formal indicial umbral calculus to Mellin--Barnes umbral transmutation theory. The starting point is the classical umbral representation of \(J_0\), which converts products of Bessel functions into Gaussian integrals involving fractional powers of sums of independent umbral operators. While this formal procedure reproduces the correct results in suitable convergence chambers, it may lead to non-admissible hypergeometric expansions at physically relevant parameter values. The cubic moment provides the basic example: the naive three-clock expansion gives an Appell \(F_4\) series whose convergence domain excludes the equilateral case, although the integral itself is finite.

	We show that this obstruction is removed by replacing the purely formal expansion with an analytic umbral transmutation. In this setting, exponential umbral pairings are interpreted through Mellin--Barnes integrals, and Ramanujan's Master Theorem acts as an inverse selection principle for the spectral ground state, or clock, associated with a given Bessel product. The factorisation \(J_0^3=J_0J_0^2\) produces two distinct clocks and reduces the cubic full-line moment to a one-dimensional Barnes integral, equivalently to a Meijer \(G\)-function. This gives the classical value of the cubic Bessel moment and clarifies why the divergent Appell realisation is only a local representation of a globally meaningful umbral identity.

	The same mechanism is then applied to scaled cubic products and to the fourth Bessel moment. For \(\int_{-\infty}^{\infty}J_0(ax)J_0(bx)^2\,\mathrm d x\), the dependence on the ratio of scales is encoded by a single effective one-clock transmutation function. For \(J_0^4\), the symmetric contraction of two square clocks yields a Meijer--Barnes representation whose right-residue expansion produces a harmonic-number series involving central binomial coefficients. The fifth moment marks the first genuinely higher-rank case: the natural umbral grouping leads to a bivariate Barnes transmutation rather than to an ordinary Meijer \(G\)-function.

	Finally, we discuss real fractional powers \(J_0^\alpha\), \(\alpha>2\), showing that the same interpretation persists beyond integer moments. In general, the corresponding Ramanujan-selected ground state is not a finite Gamma-product clock, but the moment remains a Mellin value of an analytically selected spectral state. The resulting picture identifies Bessel moments as values of effective umbral transmutations and separates the global analytic meaning of the umbral representation from the local convergence properties of its hypergeometric residue expansions.
\end{abstract}


\end{frontmatter}

%%%%%%%%%%%%%%%%%%%%%%%%%%%%%%%%%%%%%%%%%%%%%%%%%%%%%%%
\section{Introduction}\label{sec:introduction}

The interest in umbral calculus, in its distinctive indicial formulation, was initially motivated by the evaluation of infinite integrals~\cite{BabusciDattoli2011RMT}, within a formalism naturally encoding the Ramanujan's Master Theorem~\cite{Hardy1940Ramanujan}.

This program was carried out, for instance in~\cite{BabusciEtAl2012Definite,BabusciEtAl2019Methods}, using integrals involving Bessel functions as benchmark examples. The guiding paradigm is the use of the Bessel umbral image
\begin{dmath*}
	J_0(x) = {
		\exp\left(-\mathfrak u \frac{x^2}{4}\right)\varphi(0),
	\qquad
	\mathfrak u^r \varphi(0) =
		\frac{1}{\Gamma(r+1)}
	}
\end{dmath*},
where $\mathfrak u$ is the indicial umbral operator.
One can then apply ordinary Gaussian integration methods.

In the following, we apply this method to integrals involving powers of Bessel functions, analyse the limits of the method, and propose a possible way out based on the point of view recently developed in~\cite{Riccianalytic2026}.

The starting example is the evaluation of integrals containing the product \(J_0(ax)J_0(bx)\). Using independent umbral vacua, products of Bessel functions may be reduced to Gaussian integrals involving fractional powers of sums of umbral operators:
\begin{dmath*}
	J_0(ax) = {
		\exp\left(-\mathfrak u_1 \frac{a^2x^2}{4}\right)\varphi_1(0),
	\qquad
	J_0(bx) =
		\exp\left(-\mathfrak u_2 \frac{b^2x^2}{4}\right)\varphi_2(0)
		}
\end{dmath*}.
Thus, for constant \(a,b\), one finds
\begin{dmath}\label{eq:formal_I3}
	I_2(a,b) =
		\int_{-\infty}^{\infty}
			J_0(ax)J_0(bx)
		\,\mathrm dx
		=
		2\sqrt{\pi}
		\left(
			\mathfrak u_1 a^2+\mathfrak u_2 b^2
		\right)^{-1/2}
		\varphi_1(0)\varphi_2(0)
\end{dmath}.
Assuming $\abs{a}>\abs{b}$, this gives
\begin{dmath*}
	I_2(a,b) =
		\frac{2}{\abs{a}}
		\sum_{r=0}^{\infty}
			\frac{\left(\frac{1}{2}\right)_r \left(\frac{1}{2}\right)_r}{(r!)^2}
			\left(
				\frac{b^2}{a^2}
			\right)^r
\end{dmath*}.
Equivalently,
\begin{dmath*}
	I_2(a,b) =
		\frac{2}{\abs{a}}\,
		{}_2F_1
		\left(
			\frac{1}{2},\frac{1}{2};
			1;
			\frac{b^2}{a^2}
		\right)
\end{dmath*}.
Using
\begin{dmath*}
	{}_2F_1
	\left(
		\frac{1}{2},\frac{1}{2};
		1;
		k^2
	\right)
	=
	\frac{2}{\pi}K(k)
\end{dmath*},
the previous integral can be written as
\begin{dmath*}
	I_2(a,b) = {
		\frac{4}{\pi A}
		K\left(\frac{B}{A}\right),
	\qquad
	A =
		\max(\abs{a},\abs{b}),
	\qquad
	B =
		\min(\abs{a},\abs{b})
		}
\end{dmath*}.

The correctness of the result is readily verified. As a sanity check, we note that the series diverges for \(\abs{a}=\abs{b}\), coherently with the fact that the integral \(I_2(1,1)\) is not finite.

We now consider
\begin{dmath*}
	I_3(\alpha,\beta,\gamma) =
		\int_{-\infty}^{\infty}
			J_0(\alpha x)J_0(\beta x)J_0(\gamma x)
		\,\mathrm dx
\end{dmath*}.
Using three independent vacua, we can write the integral as
\begin{dmath*}
	I_3(\alpha,\beta,\gamma) =
		2\sqrt{\pi}
		\left(
			\mathfrak u_1\alpha^2
			+
			\mathfrak u_2\beta^2
			+
			\mathfrak u_3\gamma^2
		\right)^{-1/2}
		\prod_{j=1}^{3}\varphi_j(0)
\end{dmath*}.
Assuming
\begin{dmath}
	\abs{\alpha}>\abs{\beta}+\abs{\gamma}
\end{dmath},
we eventually obtain
\begin{dmath*}
	I_3(\alpha,\beta,\gamma) =
		\frac{2}{\abs{\alpha}}\,
		F_4
		\left(
			\frac{1}{2},
			\frac{1}{2};
			1,
			1;
			\frac{\beta^2}{\alpha^2},
			\frac{\gamma^2}{\alpha^2}
		\right)
\end{dmath*}.
Here \(F_4\) denotes the Appell hypergeometric function
\begin{dmath*}
	F_4(a,b;c,c';x,y)
		=
		\sum_{m=0}^{\infty}
		\sum_{n=0}^{\infty}
			\frac{
				(a)_{m+n}(b)_{m+n}
			}{
				(c)_m(c')_{n}
			}
			\frac{x^m}{m!}
			\frac{y^n}{n!}
\end{dmath*},
which converges for
\begin{dmath*}
	\sqrt{\abs{x}}+\sqrt{\abs{y}}<1
\end{dmath*},
consistently with the assumption above.

The conclusion of the formal umbral evaluation would therefore seem to be that the integral
\begin{dmath*}
	I_3(1,1,1)
		=
		\int_{-\infty}^{\infty}
			J_0(x)^3
		\,\mathrm dx
\end{dmath*}
does not provide a finite result. This, however, is not true: the integral sums to 
\begin{dmath*}
	\int_{-\infty}^{\infty}
			J_0(x)^3
		\,\mathrm d x
	= {
	\frac{\Gamma(1/6)^2}{2^{2/3}\sqrt3\,\pi^{3/2}\Gamma(5/6)}
	\simeq 1.79288\ldots
	}
\end{dmath*}
(see e.g. \cite{Watson1944Bessel, gradshteyn2014table, bailey2008elliptic}), as can also be verified with an actual numerical computation.

This is precisely where the naive expansion, deriving from the solution encoded in the formal expression above, reveals its weakness, although it does produce correct results under suitable conditions on the triple \((\alpha,\beta,\gamma)\).

In the forthcoming sections, we provide a broader perspective based on the umbral analytic framework developed in~\cite{Riccianalytic2026}, which resolves the ambiguity through Mellin--Barnes analytic continuation.

%%%%%%%%%%%%%%%%%%%%%%%%%%%%%%%%%%%%%%%%%%%%%%%%%%%%%%%
\section{The analytic umbral framework}

In the analytic umbral framework, the umbral operator $\mathfrak u$ is reinterpreted as a functional in the Borel additive complex variable u, namely:
\begin{dmath*}
	\mathfrak u \mapsto \ee^u	
\end{dmath*}.
For a suitable $f$, the formal operator function $f(\zeta\mathfrak u)$ thus corresponds to the composed functional $f(\zeta\ee^u)$. After introducing the multiplicative Borel variable $z=\ee^u$, one defines the spectral \emph{jump function} associated to $f$ as
\begin{dmath*}
	J_f(t; \zeta) \defeq {
		\mathcal M_z[f(\zeta z)](-t) =
		\int_0^\infty f(\zeta z) z^{-t-1}\,\mathrm d z
	}
\end{dmath*},
where $\mathcal M$ denotes the Mellin transform operator.
The umbral pairing with a ground state function $\varphi(t)$ is then defined by the Mellin--Barnes integral
\begin{dmath}\label{eq:umbral_analytic_pairing}
	F(\zeta) = {
		\croc{f(\zeta \ee^u), \varphi} \defeq
		\frac{1}{2\pi i} \int_{\mathcal C_t} J_f(t; \zeta) \varphi(t)\,\mathrm d t
	}
\end{dmath},
where $\mathcal C_t$ is a Barnes contour separating the poles of $J_f(t;\zeta)$ from those of $\varphi(t)$.

In particular, when $f(\zeta \ee^u) = \exp(\pm\zeta\ee^u)$ we have that
\begin{dmath*}
	J_{\mathrm{exp_\pm}}(t; \zeta) = (\mp \zeta)^t \Gamma(-t)
\end{dmath*}.

The formal umbral pairing corresponds to the series obtained by closing the contour to the right and evaluating the integral by residues. The analytic and formal evaluations coincide when $F(\zeta)$ is entire, otherwise the formal evaluation only yields the (possibly divergent) expansion of $F(\zeta)$ at the origin. 

In this sense, the analytic umbral framework provides a rigorous justification and realises a genuine extension of formal umbral calculus. We remand to \cite{Riccianalytic2026} for further details.

In typical applications, given $F(\zeta)$ and $\varphi(t)$, one aims at identifying the umbral Borel functional that makes \cref{eq:umbral_analytic_pairing} satisfied.  By reverting the point of view and fixing the jump function entering \cref{eq:umbral_analytic_pairing} to $J_{\mathrm{exp_-}}(t; \zeta)$, we can interpret the analytic pairing as a distinguished transform applied to the ground state -- typically called \emph{spectral clock} in this context --, namely the exponential umbral transmutation of $\varphi(t)$:
\begin{dmath*}
	F(\zeta) = {
		\mathfrak U_{\mathrm{exp_-}}[\varphi](\zeta) \defeq
		\frac{1}{2\pi i} \int_{\mathcal C_t} \zeta^t \Gamma(-t) \varphi(t)\,\mathrm d t
	}
\end{dmath*}.
The aim is now to identify the spectral clock that satisfies the previous identity. This problem can often be solved by relying on Ramanujan's Master Theorem. 

This is the viewpoint adopted in the following sections.

%%%%%%%%%%%%%%%%%%%%%%%%%%%%%%%%%%%%%%%%%%%%%%%%%%%%%%%
\section{Analytic umbral transmutation and Ramanujan selection}

We use exponential umbral transmutations of the form
\begin{dmath*}
	F(y)
	= {
	\mathfrak U_{\mathrm{exp_-}}[\Phi](y) =
	\left\langle
		\exp(-y\ee^u),
		\Phi(t)
	\right\rangle
	}
\end{dmath*}.
Formally, this means that
\begin{dmath*}
	F(y)
	=
	\sum_{m=0}^{\infty}
	\frac{(-y)^m}{m!}
	\Phi(m)
\end{dmath*}.
Ramanujan's Master Theorem provides the inverse analytic selection rule for the ground state (for a modern account of Ramanujan's Master Theorem and its Mellin-transform interpretation, see \cite{AmdeberhanEtAl2012RMT}). Namely, if the Mellin transform of \(F\) exists in the relevant strip, then
\begin{dmath*}
	\mathcal M[F](s)
	=
	\Gamma(s)\Phi(-s)
\end{dmath*},
where
\begin{dmath*}
	\mathcal M[F](s)
	=
	\int_0^\infty
	y^{s-1}F(y)
	\,\mathrm d y
\end{dmath*}.
Thus the Mellin-admissible continuation of the coefficient clock $\Phi(m)$ is
\begin{dmath*}
	\Phi(t)
	=
	\frac{\mathcal M[F](-t)}{\Gamma(-t)}
\end{dmath*}.
This is the sense in which Ramanujan's theorem acts as the inverse form of the analytic umbral transmutation.

For Bessel moments it is convenient to use
\begin{dmath*}
	y
	=
	x^2
\end{dmath*}.
If
\begin{dmath*}
	F_\alpha(y)
	\defeq
	J_0(\sqrt y)^\alpha
\end{dmath*},
then, since \(J_0\) is even,
\begin{dmath*}
	\int_{-\infty}^{\infty}
	J_0(x)^\alpha
	\,\mathrm d x
	=
	\mathcal M[F_\alpha]\left(\frac12\right)
\end{dmath*},
whenever the integral is convergent with the chosen branch convention. Therefore the Ramanujan-selected ground state satisfies
\begin{dmath*}
	\int_{-\infty}^{\infty}
	J_0(x)^\alpha
	\,\mathrm d x
	=
	\sqrt\pi\,
	\Phi_\alpha\left(-\frac12\right)
\end{dmath*},
with
\begin{dmath*}
	\Phi_\alpha(t)
	=
	\frac{
		\mathcal M[J_0(\sqrt y)^\alpha](-t)
	}{
		\Gamma(-t)
	}
\end{dmath*}.
For integer powers, this state may reduce to a finite Barnes or Meijer-type object. For generic fractional powers, it should instead be regarded as a genuinely Ramanujan-selected spectral state.

%%%%%%%%%%%%%%%%%%%%%%%%%%%%%%%%%%%%%%%%%%%%%%%%%%%%%%%
\section{The elementary and square Bessel clocks}

Let
\begin{dmath*}
	\mathcal C_0(z)
	= {
	J_0(2\sqrt z)
	=
	\sum_{m=0}^{\infty}
	\frac{(-z)^m}{\Gamma(1+m)^2}
	}
\end{dmath*}.
Thus
\begin{dmath*}
	J_0(x)
	=
	\mathcal C_0\left(\frac{x^2}{4}\right)
\end{dmath*}.
We use standard notation and normalisations for Bessel functions as in \cite{Watson1944Bessel,OlverEtAl2010NIST,NIST:DLMF}.

The elementary Bessel--Clifford clock is
\begin{dmath*}
	\varphi_1(t)
	=
	\frac{1}{\Gamma(1+t)}
\end{dmath*},
and therefore
\begin{dmath*}
	J_0(x)
	=
	\left\langle
		\exp\left(-\frac{x^2}{4}\exp u\right),
		\varphi_1(t)
	\right\rangle
\end{dmath*}.

The square has a different clock. Since
\begin{dmath*}
	\mathcal C_0(z)^2
	=
	\sum_{m=0}^{\infty}
	(-z)^m
	\sum_{p+q=m}
	\frac{1}{\Gamma(1+p)^2\Gamma(1+q)^2}
\end{dmath*}
and
\begin{dmath*}
	\sum_{p=0}^{m}
	\binom{m}{p}^2
	=
	\binom{2m}{m}
\end{dmath*},
one obtains
\begin{dmath*}
	\mathcal C_0(z)^2
	=
	\sum_{m=0}^{\infty}
	(-z)^m
	\frac{\Gamma(2m+1)}{\Gamma(m+1)^4}
\end{dmath*}.
Writing this as
\begin{dmath*}
	\mathcal C_0(z)^2
	=
	\left\langle
		\exp(-z\exp v),
		\varphi_2(s)
	\right\rangle
\end{dmath*}
forces
\begin{dmath*}
	\frac{\varphi_2(m)}{\Gamma(1+m)}
	=
	\frac{\Gamma(2m+1)}{\Gamma(m+1)^4}
\end{dmath*}.
Therefore
\begin{dmath*}
	\varphi_2(s)
	=
	\frac{\Gamma(2s+1)}{\Gamma(s+1)^3}
\end{dmath*},
and
\begin{dmath*}
	J_0(x)^2
	=
	\left\langle
		\exp\left(-\frac{x^2}{4}\exp v\right),
		\varphi_2(s)
	\right\rangle
\end{dmath*}.

%%%%%%%%%%%%%%%%%%%%%%%%%%%%%%%%%%%%%%%%%%%%%%%%%%%%%%%
\section{The cubic moment as a two-clock contraction}

We now consider
\begin{dmath*}
	I_3
	\defeq
	\int_{-\infty}^{\infty}
	J_0(x)^3
	\,\mathrm d x
\end{dmath*}.
The useful factorisation is
\begin{dmath*}
	J_0(x)^3
	=
	J_0(x)J_0(x)^2
\end{dmath*}.
Using the clocks \(\varphi_1\) and \(\varphi_2\), this gives
\begin{dmath*}
	J_0(x)^3
	=
	\left\langle
		\left\langle
			\exp\left(-\frac{x^2}{4}(\exp u+\exp v)\right),
			\varphi_1(t)\,
			\varphi_2(s)
		\right\rangle
	\right\rangle
\end{dmath*}.
After the usual admissible exchange between the spatial integral and the analytic pairing, the Gaussian integral gives
\begin{dmath*}
	\int_{-\infty}^{\infty}
	\exp\left(-\frac{x^2}{4}(\exp u+\exp v)\right)
	\,\mathrm d x
	=
	2\sqrt\pi\,(\exp u+\exp v)^{-1/2}
\end{dmath*}.
Hence
\begin{dmath*}
	I_3
	=
	2\sqrt\pi
	\left\langle
		\left\langle
			(\exp u+\exp v)^{-1/2},
			\varphi_1(t)\,
			\varphi_2(s)
		\right\rangle
	\right\rangle
\end{dmath*}.
This formula is analogous to the formal pairing of \cref{eq:formal_I3} for $a=b=1$, but its analytic character enables to bypass the limitations of the formal approach.

Specifically, we use the Mellin--Barnes decomposition
\begin{dmath*}
	(a+b)^{-1/2}
	=
	\frac{1}{\sqrt\pi}
	\frac{1}{2\pi i}
	\int_{\mathcal C_r}
	\Gamma(-r)
	\Gamma\left(r+\frac12\right)
	a^r b^{-r-1/2}
	\,\mathrm d r
\end{dmath*},
where the contour separates the poles of \(\Gamma(-r)\) from those of \(\Gamma(r+1/2)\). With \(a=\exp u\) and \(b=\exp v\), the two clocks evaluate at complementary spectral values:
\begin{dmath*}
	\left\langle
		\exp(ru),
		\varphi_1(t)
	\right\rangle
	=
	\varphi_1(r)
\end{dmath*},
and
\begin{dmath*}
	\left\langle
		\exp\left(\left(-r-\frac12\right)v\right),
		\varphi_2(s)
	\right\rangle
	=
	\varphi_2\left(-r-\frac12\right)
\end{dmath*}.
Thus
\begin{dmath*}
	I_3
	=
	\frac{2}{2\pi i}
	\int_{\mathcal C_r}
	\Gamma(-r)
	\Gamma\left(r+\frac12\right)
	\varphi_1(r)
	\varphi_2\left(-r-\frac12\right)
	\,\mathrm d r
\end{dmath*}.
Since
\begin{dmath*}
	\varphi_1(r)
	=
	\frac{1}{\Gamma(1+r)}
\end{dmath*}
and
\begin{dmath*}
	\varphi_2\left(-r-\frac12\right)
	=
	\frac{\Gamma(-2r)}{\Gamma\left(\frac12-r\right)^3}
\end{dmath*},
we obtain
\begin{dmath*}
	I_3
	=
	\frac{2}{2\pi i}
	\int_{\mathcal C_r}
	\Gamma(-r)
	\Gamma\left(r+\frac12\right)
	\frac{1}{\Gamma(1+r)}
	\frac{\Gamma(-2r)}{\Gamma\left(\frac12-r\right)^3}
	\,\mathrm d r
\end{dmath*}.
Using
\begin{dmath*}
	\Gamma(-2r)
	=
	2^{-2r-1}\pi^{-1/2}
	\Gamma(-r)
	\Gamma\left(\frac12-r\right)
\end{dmath*},
this becomes
\begin{dmath*}
	I_3
	=
	\frac{1}{\sqrt\pi}
	\frac{1}{2\pi i}
	\int_{\mathcal C_r}
	4^{-r}
	\frac{
		\Gamma(-r)^2
		\Gamma\left(r+\frac12\right)
	}{
		\Gamma(1+r)
		\Gamma\left(\frac12-r\right)^2
	}
	\,\mathrm d r
\end{dmath*}.
Equivalently, after setting \(r=-t\),
\begin{dmath*}
	I_3
	=
	\frac{1}{\sqrt\pi}
	\frac{1}{2\pi i}
	\int_{\mathcal C_t}
	4^t
	\frac{
		\Gamma(t)^2
		\Gamma\left(\frac12-t\right)
	}{
		\Gamma(1-t)
		\Gamma\left(\frac12+t\right)^2
	}
	\,\mathrm d t
\end{dmath*}.

This Barnes integral can also be directly obtained from the Mellin--Parseval formula. Indeed,
\begin{dmath*}
	\mathcal M[J_0](s)
	=
	2^{s-1}
	\frac{\Gamma(s/2)}{\Gamma(1-s/2)}
\end{dmath*},
and
\begin{dmath*}
	\mathcal M[J_0^2](s)
	=
	\frac{1}{2\sqrt\pi}
	\frac{
		\Gamma(s/2)
		\Gamma(1/2-s/2)
	}{
		\Gamma(1-s/2)^2
	}
\end{dmath*}.
Therefore
\begin{dmath*}
	I_3
	=
	2
	\frac{1}{2\pi i}
	\int_{\mathcal C_s}
	\mathcal M[J_0](s)
	\mathcal M[J_0^2](1-s)
	\,\mathrm d s
\end{dmath*},
and the change of variable \(s=2t\) gives precisely the preceding Barnes integral.

In Meijer--Barnes notation,
\begin{dmath*}
	I_3
	=
	\frac{1}{\sqrt\pi}
	G_{3,3}^{2,1}
	\left(
		\frac14
		\left|
		\begin{matrix}
			\frac12,\frac12,\frac12\\
			0,0,0
		\end{matrix}
		\right.
	\right)
\end{dmath*}.
Our conventions for Meijer--Barnes integrals and Meijer \(G\)-functions are those of \cite{PrudnikovBrychkovMarichev1986Vol2,NIST:DLMF}; for a broader Mellin--Barnes framework, see also \cite{MathaiSaxena1978HFunction}.

The numeric value of $I_3$, as verified with the aid of Mathematica{\texttrademark}, is
\begin{dmath*}
	\frac{1}{\sqrt\pi}
	G_{3,3}^{2,1}
	\left(
		\frac14
		\left|
		\begin{matrix}
			\frac12,\frac12,\frac12\\
			0,0,0
		\end{matrix}
		\right.
	\right)
	\simeq 1.79288\ldots
\end{dmath*}.
This numerically agrees with classical Watson's evaluation reported in \cref{sec:introduction}:
\begin{dmath*}
	I_3
	=
	\frac{\Gamma(1/6)^2}{2^{2/3}\sqrt3\,\pi^{3/2}\Gamma(5/6)}
\end{dmath*}.

\begin{remark}
	The Mellin--Barnes value of the cubic Bessel moment can be written in the hypergeometric form
	\begin{dmath*}
		I_3 =
			\sqrt{3}\,
			{}_3F_2
			\left(
				\begin{matrix}
					\frac12,\frac12,\frac12\\
					1,1
				\end{matrix}
				;
				\frac14
			\right)
	\end{dmath*}.
	By Clausen's identity \cite{Bailey1935,AndrewsAskeyRoy1999} ,
	\begin{dmath*}
		{}_3F_2
		\left(
			\begin{matrix}
				\frac12,\frac12,\frac12\\
				1,1
			\end{matrix}
			;
			z
		\right)
		=
		\left[
			{}_2F_1
			\left(
				\begin{matrix}
					\frac14,\frac14\\
					1
				\end{matrix}
				;
				z
			\right)
		\right]^2
	\end{dmath*}.
	Therefore
	\begin{dmath*}
		I_3 =
			\sqrt{3}
			\left[
				{}_2F_1
				\left(
					\begin{matrix}
						\frac14,\frac14\\
						1
					\end{matrix}
					;
					\frac14
				\right)
			\right]^2
	\end{dmath*}.
	The singular-modulus evaluation \cite{BorweinBorwein1987}
	\begin{dmath*}
		{}_2F_1
		\left(
			\begin{matrix}
				\frac14,\frac14\\
				1
			\end{matrix}
			;
			\frac14
		\right)
		=
		\frac{
			\Gamma\left(\frac16\right)
		}{
			2^{1/3}\sqrt{3}\,\pi^{3/4}
			\Gamma\left(\frac56\right)^{1/2}
		}
	\end{dmath*}
	then gives
	\begin{dmath*}
		I_3 =
			\frac{
				\Gamma\left(\frac16\right)^2
			}{
				2^{2/3}\sqrt{3}\,\pi^{3/2}
				\Gamma\left(\frac56\right)
			}
	\end{dmath*}.
	This identifies the Mellin--Barnes continuation of the formal cubic umbral expression with the classical closed value of the three-Bessel moment.
\end{remark}

%%%%%%%%%%%%%%%%%%%%%%%%%%%%%%%%%%%%%%%%%%%%%%%%%%%%%%%
\section{Scaled cubic products and an effective one-clock function}

The same computation applies to the scaled product integral
\begin{dmath*}
	I(a,b)
	\defeq
	\int_{-\infty}^{\infty}
	J_0(ax)J_0(bx)^2
	\,\mathrm d x
\end{dmath*}.
The two-clock representation is
\begin{dmath*}
	J_0(ax)J_0(bx)^2
	=
	\left\langle
		\left\langle
			\exp\left(
				-\frac{x^2}{4}
				\left(
					a^2\exp u
					+
					b^2\exp v
				\right)
			\right),
			\varphi_1(t)\,
			\varphi_2(s)
		\right\rangle
	\right\rangle
\end{dmath*}.
The Gaussian integral and the same Barnes decomposition yield
\begin{dmath*}
	I(a,b)
	=
	\frac{1}{\sqrt\pi\,b}
	\frac{1}{2\pi i}
	\int_{\mathcal C_r}
	\left(
		\frac{a^2}{4b^2}
	\right)^r
	\frac{
		\Gamma(-r)^2
		\Gamma\left(\frac12+r\right)
	}{
		\Gamma(1+r)
		\Gamma\left(\frac12-r\right)^2
	}
	\,\mathrm d r
\end{dmath*}.
Define
\begin{dmath*}
	\mathcal H(\zeta)
	\defeq
	\frac{1}{\sqrt\pi}
	\frac{1}{2\pi i}
	\int_{\mathcal C_r}
	\zeta^r
	\frac{
		\Gamma(-r)^2
		\Gamma\left(\frac12+r\right)
	}{
		\Gamma(1+r)
		\Gamma\left(\frac12-r\right)^2
	}
	\,\mathrm d r
\end{dmath*}.
Then
\begin{dmath*}
	I(a,b)
	=
	\frac{1}{b}
	\mathcal H\left(\frac{a^2}{4b^2}\right)
\end{dmath*}.
Equivalently,
\begin{dmath*}
	\mathcal H(\zeta)
	=
	\frac{1}{\sqrt\pi}
	G_{3,3}^{2,1}
	\left(
		\zeta
		\left|
		\begin{matrix}
			\frac12,\frac12,\frac12\\
			0,0,0
		\end{matrix}
		\right.
	\right)
\end{dmath*},
so that
\begin{dmath*}
	I(a,b)
	=
	\frac{1}{\sqrt\pi\,b}
	G_{3,3}^{2,1}
	\left(
		\frac{a^2}{4b^2}
		\left|
		\begin{matrix}
			\frac12,\frac12,\frac12\\
			0,0,0
		\end{matrix}
		\right.
	\right)
\end{dmath*},
with the usual Barnes-contour prescription and analytic continuation in the ratio \(a/b\) where necessary.

The same function can be interpreted as an effective one-clock transmutation:
\begin{dmath*}
	\mathcal H(\zeta)
	=
	\left\langle
		\exp(-\zeta\exp u),
		\varphi_{\mathrm{eff}}(t)
	\right\rangle
\end{dmath*},
where
\begin{dmath*}
	\varphi_{\mathrm{eff}}(t)
	=
	\frac{1}{\sqrt\pi}
	\frac{
		\Gamma(-t)
		\Gamma\left(\frac12+t\right)
	}{
		\Gamma(1+t)
		\Gamma\left(\frac12-t\right)^2
	}
\end{dmath*}.
This effective clock contains the spectral convolution of the elementary clock \(\varphi_1\) and the square clock \(\varphi_2\).

%%%%%%%%%%%%%%%%%%%%%%%%%%%%%%%%%%%%%%%%%%%%%%%%%%%%%%%
\section{The fourth moment}

We now consider
\begin{dmath*}
	I_4
	\defeq
	\int_{-\infty}^{\infty}
	J_0(x)^4
	\,\mathrm d x
\end{dmath*}.
The natural grouping is symmetric:
\begin{dmath*}
	J_0(x)^4
	=
	J_0(x)^2J_0(x)^2
\end{dmath*}.
Using the square clock \(\varphi_2\), we obtain
\begin{dmath*}
	J_0(x)^4
	=
	\left\langle
		\left\langle
			\exp\left(
				-\frac{x^2}{4}(\exp u+\exp v)
			\right),
			\varphi_2(t)\,
			\varphi_2(s)
		\right\rangle
	\right\rangle
\end{dmath*}.
The same Gaussian integration and Barnes decomposition give
\begin{dmath*}
	I_4
	=
	2
	\frac{1}{2\pi i}
	\int_{\mathcal C_r}
	\Gamma(-r)
	\Gamma\left(r+\frac12\right)
	\varphi_2(r)
	\varphi_2\left(-r-\frac12\right)
	\,\mathrm d r
\end{dmath*}.
Substitution of
\begin{dmath*}
	\varphi_2(r)
	=
	\frac{\Gamma(2r+1)}{\Gamma(r+1)^3}
\end{dmath*}
and
\begin{dmath*}
	\varphi_2\left(-r-\frac12\right)
	=
	\frac{\Gamma(-2r)}{\Gamma\left(\frac12-r\right)^3}
\end{dmath*}
yields
\begin{dmath*}
	I_4
	=
	\frac{2}{2\pi i}
	\int_{\mathcal C_r}
	\Gamma(-r)
	\Gamma\left(r+\frac12\right)
	\frac{\Gamma(2r+1)}{\Gamma(r+1)^3}
	\frac{\Gamma(-2r)}{\Gamma\left(\frac12-r\right)^3}
	\,\mathrm d r
\end{dmath*}.
Using
\begin{dmath*}
	\Gamma(2r+1)
	=
	2^{2r}\pi^{-1/2}
	\Gamma\left(r+\frac12\right)
	\Gamma(r+1)
\end{dmath*}
and
\begin{dmath*}
	\Gamma(-2r)
	=
	2^{-2r-1}\pi^{-1/2}
	\Gamma(-r)
	\Gamma\left(\frac12-r\right)
\end{dmath*},
the powers of \(2\) cancel the external factor \(2\), and one obtains
\begin{dmath*}
	I_4
	=
	\frac{1}{\pi}
	\frac{1}{2\pi i}
	\int_{\mathcal C_r}
	\frac{
		\Gamma(-r)^2
		\Gamma\left(r+\frac12\right)^2
	}{
		\Gamma(r+1)^2
		\Gamma\left(\frac12-r\right)^2
	}
	\,\mathrm d r
\end{dmath*}.

Define the effective fourth-moment function
\begin{dmath*}
	\mathcal H_4(\zeta)
	\defeq
	\frac{1}{\pi}
	\frac{1}{2\pi i}
	\int_{\mathcal C_r}
	\zeta^r
	\frac{
		\Gamma(-r)^2
		\Gamma\left(r+\frac12\right)^2
	}{
		\Gamma(r+1)^2
		\Gamma\left(\frac12-r\right)^2
	}
	\,\mathrm d r
\end{dmath*}.
Then
\begin{dmath*}
	I_4
	=
	\mathcal H_4(1)
\end{dmath*}.
As a one-clock transmutation,
\begin{dmath*}
	\mathcal H_4(\zeta)
	=
	\left\langle
		\exp(-\zeta\exp u),
		\varphi_{\mathrm{eff},4}(t)
	\right\rangle
\end{dmath*},
where
\begin{dmath*}
	\varphi_{\mathrm{eff},4}(t)
	=
	\frac{1}{\pi}
	\frac{
		\Gamma(-t)
		\Gamma\left(t+\frac12\right)^2
	}{
		\Gamma(t+1)^2
		\Gamma\left(\frac12-t\right)^2
	}
\end{dmath*}.
In Meijer--Barnes notation,
\begin{dmath*}
	\mathcal H_4(\zeta)
	=
	\frac{1}{\pi}
	G_{4,4}^{2,2}
	\left(
		\zeta
		\left|
		\begin{matrix}
			\frac12,\frac12,\frac12,\frac12\\
			0,0,0,0
		\end{matrix}
		\right.
	\right)
\end{dmath*}.
Therefore
\begin{dmath*}
	\int_{-\infty}^{\infty}
	J_0(x)^4
	\,\mathrm d x
	=
	\frac{1}{\pi}
	G_{4,4}^{2,2}
	\left(
		1
		\left|
		\begin{matrix}
			\frac12,\frac12,\frac12,\frac12\\
			0,0,0,0
		\end{matrix}
		\right.
	\right)
\end{dmath*},
and numerically
\begin{dmath*}
	\int_{-\infty}^{\infty}
	J_0(x)^4
	\,\mathrm d x
	=
	1.80545\ldots
\end{dmath*}.

\proofstep{Residue expansion.}
Keep \(\zeta\) as an auxiliary parameter and assume first that \(|\zeta|<1\). Closing the Barnes contour to the right, the enclosed poles are the double poles at \(r=n\), \(n=0,1,2,\dots\), coming from \(\Gamma(-r)^2\). With the usual right-closure orientation one obtains
\begin{dmath*}
	\mathcal H_4(\zeta)
	=
	-\frac{1}{\pi}
	\sum_{n=0}^{\infty}
	\operatorname*{Res}_{r=n}
	\left[
		\zeta^r
		\frac{
			\Gamma(-r)^2
			\Gamma\left(r+\frac12\right)^2
		}{
			\Gamma(r+1)^2
			\Gamma\left(\frac12-r\right)^2
		}
	\right]
\end{dmath*}.
The residue calculation is a standard double-pole Barnes computation; compare the Meijer \(G\)-function residue conventions in \cite{NIST:DLMF, PrudnikovBrychkovMarichev1986Vol2}.

A direct double-pole expansion gives
\begin{dmath*}
	\operatorname*{Res}_{r=n}
	\left[
		\zeta^r
		\frac{
			\Gamma(-r)^2
			\Gamma\left(r+\frac12\right)^2
		}{
			\Gamma(r+1)^2
			\Gamma\left(\frac12-r\right)^2
		}
	\right]
	=
	\frac{\binom{2n}{n}^4}{256^n}
	\zeta^n
	\left[
		\log\zeta
		+
		8\left(
			H_{2n}-H_n-\log2
		\right)
	\right]
\end{dmath*}.
Therefore
\begin{dmath*}
	\mathcal H_4(\zeta)
	=
	-\frac{1}{\pi}
	\sum_{n=0}^{\infty}
	\frac{\binom{2n}{n}^4}{256^n}
	\zeta^n
	\left[
		\log\zeta
		+
		8\left(
			H_{2n}-H_n-\log2
		\right)
	\right]
\end{dmath*}.
Taking the limit \(\zeta\to1^{-}\), the logarithmic term disappears and
\begin{dmath*}
	I_4
	=
	\mathcal H_4(1)
	=
	\frac{8}{\pi}
	\sum_{n=0}^{\infty}
	\frac{\binom{2n}{n}^4}{256^n}
	\left(
		\log2+H_n-H_{2n}
	\right)
\end{dmath*}.
The convergence follows from
\begin{dmath*}
	\frac{\binom{2n}{n}^4}{256^n}
	\sim
	\frac{1}{\pi^2n^2}
\end{dmath*}
and
\begin{dmath*}
	\log2+H_n-H_{2n}
	\sim
	\frac{1}{4n}
\end{dmath*},
so the summand is \(O(n^{-3})\). Numerically, this residue expansion gives
\begin{dmath*}
	\frac{8}{\pi}
	\sum_{n=0}^{\infty}
	\frac{\binom{2n}{n}^4}{256^n}
	\left(
		\log2+H_n-H_{2n}
	\right)
	=
	1.80545\ldots
\end{dmath*}.

%%%%%%%%%%%%%%%%%%%%%%%%%%%%%%%%%%%%%%%%%%%%%%%%%%%%%%%
\section{Beyond Meijer \(G\): the fifth Bessel moment as a bivariate Barnes transmutation}\label{sec:fifth-bessel-bivariate}

	The preceding examples show that the cubic and quartic Bessel moments admit particularly efficient one-dimensional Barnes representations. This is not because the analytic umbral method is intrinsically restricted to third and fourth powers, but because in those cases the relevant product can be reduced to a binary fusion of two single-state analytic transmutations. The binomial Barnes representation then introduces only one auxiliary Barnes variable, and the result is naturally expressible as a Meijer \(G\)-function.

	The fifth moment is the first case where the same mechanism naturally produces a genuinely multivariate Barnes object. We group
	\begin{dmath*}
		J_0(x)^5
		=
		J_0(x)^2\,J_0(x)^2\,J_0(x)
	\end{dmath*}
	and use two copies of the square clock
	\begin{dmath*}
		\varphi_2(t)
		=
		\frac{\Gamma(2t+1)}{\Gamma(t+1)^3}
	\end{dmath*}
	together with the elementary clock
	\begin{dmath*}
		\varphi_1(t)
		=
		\frac{1}{\Gamma(1+t)}
	\end{dmath*}.
	In umbral form, this means that
	\begin{dmath*}
		J_0(x)^2
		=
		\left\langle
			\exp\left(
				-\frac{x^2}{4}\exp u
			\right),
			\varphi_2(t)
		\right\rangle
	\end{dmath*}
	and
	\begin{dmath*}
		J_0(x)
		=
		\left\langle
			\exp\left(
				-\frac{x^2}{4}\exp u
			\right),
			\varphi_1(t)
		\right\rangle
	\end{dmath*}.
	Therefore the fifth power is represented by a three-clock pairing. After the Gaussian integration in the spatial variable, the kernel becomes
	\begin{dmath*}
		\left(
			\exp u
			+
			\exp v
			+
			\exp w
		\right)^{-1/2}
	\end{dmath*}.
	This is the point at which the fifth moment departs from the Meijer \(G\)-type behaviour of the cubic and quartic cases.

	The appropriate Barnes representation is now the multinomial one:
	\begin{dmath*}
		(a+b+c)^{-1/2}
		=
		\frac{1}{\Gamma(1/2)}
		\frac{1}{(2\pi i)^2}
		\int_{\mathcal C_r}
		\int_{\mathcal C_s}
		\Gamma(-r)
		\Gamma(-s)
		\Gamma\left(
			\frac12+r+s
		\right)
		a^r b^s c^{-1/2-r-s}
		\,\mathrm dr\,\mathrm ds
	\end{dmath*}.
	Choosing
	\begin{dmath*}
		a = {
			\exp u,
			\qquad
			b =
			\exp v,
			\qquad
			c =
			\exp w
		}
	\end{dmath*}
	shows that the three clocks are evaluated respectively at
	\begin{dmath*}
		r,
		\qquad
		s,
		\qquad
		-\frac12-r-s
	\end{dmath*}.
	Consequently the full fifth moment can be written formally as
	\begin{dmath*}
		I_5
		=
		2
		\frac{1}{(2\pi i)^2}
		\int_{\mathcal C_r}
		\int_{\mathcal C_s}
		\Gamma(-r)
		\Gamma(-s)
		\Gamma\left(
			\frac12+r+s
		\right)
		\varphi_2(r)
		\varphi_2(s)
		\varphi_1\left(
			-\frac12-r-s
		\right)
		\,\mathrm dr\,\mathrm ds
	\end{dmath*}.
	Substituting the two clocks gives the cleaner bivariate Barnes representation
	\begin{dmath*}
		I_5
		=
		2
		\frac{1}{(2\pi i)^2}
		\int_{\mathcal C_r}
		\int_{\mathcal C_s}
		\Gamma(-r)
		\Gamma(-s)
		\frac{
			\Gamma\left(
				\frac12+r+s
			\right)
		}{
			\Gamma\left(
				\frac12-r-s
			\right)
		}
		\frac{\Gamma(2r+1)}{\Gamma(r+1)^3}
		\frac{\Gamma(2s+1)}{\Gamma(s+1)^3}
		\,\mathrm dr\,\mathrm ds
	\end{dmath*}.
	This formula should be read with the same Barnes-contour conventions used in the preceding one-dimensional cases, but now with two independent contour variables.

	It is useful to introduce two auxiliary parameters and regard the last expression as the value at \((1,1)\) of a bivariate Barnes function. Define
	\begin{dmath*}
		\mathcal H_5(\zeta,\eta)
		=
		\frac{1}{(2\pi i)^2}
		\int_{\mathcal C_r}
		\int_{\mathcal C_s}
		\zeta^r
		\eta^s
		\Gamma(-r)
		\Gamma(-s)
		\Phi_{5}^{\mathrm{eff}}(r,s)
		\,\mathrm dr\,\mathrm ds =
		=
		\left\langle \left\langle \exp(\zeta \ee^u + \eta \ee^v), \Phi_{5}^{\mathrm{eff}} \right\rangle \right\rangle
	\end{dmath*},
	where the bivariate effective ground state is
	\begin{dmath*}
		\Phi_{5}^{\mathrm{eff}}(r,s)
		=
		2\,
		\frac{
			\Gamma\left(
				\frac12+r+s
			\right)
		}{
			\Gamma\left(
				\frac12-r-s
			\right)
		}
		\frac{\Gamma(2r+1)}{\Gamma(r+1)^3}
		\frac{\Gamma(2s+1)}{\Gamma(s+1)^3}
	\end{dmath*}.
	Then
	\begin{dmath*}
		I_5
		=
		\mathcal H_5(1,1)
	\end{dmath*}.
	In this sense the fifth moment is still governed by an effective umbral state, but the state is now bivariate. The ordinary Meijer \(G\)-function is replaced by a two-variable Mellin--Barnes function.

	Closing both Barnes contours to the right gives local Horn-type expansions for \(\mathcal H_5(\zeta,\eta)\), rather than a directly convergent formula at the physical point \((1,1)\). Indeed, the right poles are
	\begin{dmath*}
		r = {
		m,
		\qquad
		s=n,
		\qquad
		m,n\in\mathbb N
		}
	\end{dmath*}
	and the resulting residue expansion is a local expansion in the variables \(\zeta\) and \(\eta\). Unlike the fourth-moment case, the endpoint \((\zeta,\eta)=(1,1)\) does not lie in the domain where this naive right-residue expansion gives an efficient numerical summation. Thus the formula
	\begin{dmath*}
		I_5
		=
		\mathcal H_5(1,1)
	\end{dmath*}
	should be interpreted as an analytic-continuation statement for a bivariate Barnes transmutation, not as a directly useful residue series.

	This gives the natural hierarchy behind the preceding computations:
	\begin{dmath*}
		\begin{aligned}
			n=3,4
			&\colon
			\text{binary fusion}
			\longrightarrow
			\text{one effective clock}
			\longrightarrow
			\text{one-variable Meijer }G,
			\\
			n=5
			&\colon
			\text{ternary fusion}
			\longrightarrow
			\text{bivariate effective clock}
			\longrightarrow
			\text{two-variable Barnes function}.
		\end{aligned}
	\end{dmath*}
	Thus the obstruction at the fifth moment is not a failure of the analytic umbral method. Rather, it is the first appearance of the higher-rank Barnes structure which is suppressed in the cubic and quartic cases by the availability of a binary reduction.

	For actual computation, however, the bivariate Barnes representation is not the most efficient tool. One may instead use other methods, such as the angular product formula for \(J_0^2\), which reduces the problem to a two-dimensional integral involving the reciprocal area of a Euclidean triangle \cite{Watson1944Bessel}.
	
	In this sense the fifth moment separates the structural content of the analytic-umbral method from the most efficient numerical representation.

%%%%%%%%%%%%%%%%%%%%%%%%%%%%%%%%%%%%%%%%%%%%%%%%%%%%%%%
\section{Fractional powers}

We now return to
\begin{dmath*}
	I_\alpha
	\defeq
	\int_{-\infty}^{\infty}
	J_0(x)^\alpha
	\,\mathrm d x
\end{dmath*},
with the principal branch convention
\begin{dmath*}
	J_0(x)^\alpha
	\defeq
	\exp\left(\alpha\,\mathrm{Log}\,J_0(x)\right)
\end{dmath*}.
For non-integer \(\alpha\), the zeros of \(J_0\) become branch points of the integrand. If \(j_0=0<j_1<j_2<\cdots\) are the non-negative zeros of \(J_0\), then for real \(\alpha\) the Mellin transform decomposes as
\begin{dmath*}
	\mathcal M[F_\alpha](s)
	=
	2
	\sum_{k=0}^{\infty}
	\ee^{i\pi\alpha k}
	\int_{j_k}^{j_{k+1}}
	x^{2s-1}
	|J_0(x)|^\alpha
	\,\mathrm d x
\end{dmath*},
because the sign of \(J_0\) alternates from interval to interval.

The convergence strip follows from the endpoint behaviour. At the origin,
\begin{dmath*}
	J_0(x)
	=
	1+O(x^2)
	\condition*{x\to0}
\end{dmath*},
so one needs
\begin{dmath*}
	\operatorname{Re}s>0
\end{dmath*}.
At infinity,
\begin{dmath*}
	J_0(x)
	\sim
	\left(\frac{2}{\pi x}\right)^{1/2}
	\cos\left(x-\frac{\pi}{4}\right)
\end{dmath*}.
For real non-integer \(\alpha\), there is no systematic cancellation sufficient to improve the absolute estimate, and the tail is absolutely integrable when
\begin{dmath*}
	2\operatorname{Re}s-1-\frac{\alpha}{2}
	<
	-1
\end{dmath*},
that is,
\begin{dmath*}
	\operatorname{Re}s
	<
	\frac{\alpha}{4}
\end{dmath*}.
Thus a natural absolute-convergence strip is
\begin{dmath*}
	0
	< {
	\operatorname{Re}s
	<
	\frac{\alpha}{4}
	}
\end{dmath*}.
In particular, the value \(s=1/2\) belongs to the strip exactly when
\begin{dmath*}
	\alpha>2
\end{dmath*}.
The zeros of \(J_0\) impose no stronger condition in this range: the zeros are simple, and locally \(|J_0(x)|^\alpha\) is integrable whenever \(\alpha>-1\).

Therefore, for real \(\alpha>2\), the principal-branch integral is absolutely convergent and
\begin{dmath*}
	I_\alpha
	=
	\sqrt\pi\,
	\Phi_\alpha\left(-\frac12\right)
\end{dmath*},
where
\begin{dmath*}
	\Phi_\alpha(t)
	=
	\frac{
		\mathcal M[J_0(\sqrt y)^\alpha](-t)
	}{
		\Gamma(-t)
	}
\end{dmath*}.
This formula is exact, but it is a structural reduction rather than an automatic closed-form evaluation. The integer cases are special because the branch obstruction disappears and the effective state may collapse to a finite Barnes kernel.

\begin{remark}
	For \(\alpha=1\), the full-line integral is conditionally convergent and equals \(2\). For \(\alpha=2\), the integral diverges logarithmically, consistently with the boundary \(\alpha>2\). The cases \(\alpha=3\) and \(\alpha=4\) are the first absolutely convergent integer cases and are treated explicitly in this article.
\end{remark}

%%%%%%%%%%%%%%%%%%%%%%%%%%%%%%%%%%%%%%%%%%%%%%%%%%%%%%%
\section{Conclusion}

We have shown how a consistent analytic reformulation of the umbral framework enables to bypass the limitations of the formal umbral approach in the determination of the Bessel moments, by reproducing and reinterpreting classical results.

 The main conclusion we can draw from the carried out analysis is structural. Products of Bessel functions can be reorganised into exponential umbral transmutations with distinct spectral clocks. Spatial integration then produces elementary algebraic kernels whose Barnes decompositions fuse the clocks into effective umbral states.

For \(J_0^3\), the relevant clocks are \(\varphi_1\) and \(\varphi_2\), and the binomial Barnes decomposition produces the one-variable Meijer--Barnes function \(\mathcal H\). For the scaled product \(J_0(ax)J_0(bx)^2\), the entire dependence on the ratio of scales enters through the single argument $\frac{a^2}{4b^2}$.
Thus the scaled moment is controlled by the same effective function \(\mathcal H\).

For \(J_0^4\), the symmetric pair of square clocks \(\varphi_2,\varphi_2\) again leads to a binary fusion. The resulting effective function \(\mathcal H_4\) is a one-variable Meijer--Barnes function, and its value at \(1\) gives the fourth Bessel moment. In this case closing the Barnes contour to the right gives a useful Abel endpoint expansion: the double poles produce a rapidly convergent harmonic-number series for \(I_4\).

The fifth moment clarifies the scope of this Meijer \(G\)-phenomenon. Grouping
\begin{dmath*}
	J_0(x)^5
	=
	J_0(x)^2\,J_0(x)^2\,J_0(x)
\end{dmath*}
still gives a perfectly natural analytic umbral representation, but the spatial integral now produces a ternary kernel. Its multinomial Barnes decomposition introduces two independent Barnes variables. Consequently the effective object is no longer a one-clock state and no longer an ordinary Meijer \(G\)-function. Instead, \(I_5\) is realised as the value
\begin{dmath*}
	I_5
	=
	\mathcal H_5(1,1)
\end{dmath*}
of a bivariate Barnes transmutation. Closing the two Barnes contours gives local Horn-type expansions rather than a directly convergent residue series at the physical point. Thus the fifth moment separates the structural analytic umbral representation from the most efficient numerical representation, which is better obtained from the angular product formula for \(J_0^2\).

The fractional-power discussion explains why the same language persists beyond integer moments. For real \(\alpha>2\), the principal-branch moment is given by the Mellin value
\begin{dmath*}
	\mathcal M
	\left[
		J_0(\sqrt y)^\alpha
	\right]
	\left(
		\frac12
	\right)
\end{dmath*}
or, equivalently, by
\begin{dmath*}
	\sqrt{\pi}\,
	\Phi_\alpha\left(
		-\frac12
	\right)
\end{dmath*},
where \(\Phi_\alpha\) is selected by Ramanujan's Master Theorem. In general this selected state is not a finite Gamma product. Nevertheless, the analytic umbral interpretation remains the same: Bessel moments are Mellin values of Ramanujan-selected ground states, and the passage from \(n=3,4\) to \(n=5\) displays the transition from rank-one Meijer--Barnes transmutations to higher-rank Barnes structures.
	
As a final remark, we observe that although the naive Appell series associated with the formal evaluation of \(I_3(1,1,1)\), namely
\begin{dmath*}
	F_4\left(
		\frac12,\frac12;
		1,1;
		1,1
	\right)
	=
	\sum_{m=0}^{\infty}
	\sum_{n=0}^{\infty}
		\frac{
			\left(\frac12\right)_{m+n}^2
		}{
			(m!)^2(n!)^2
		}
\end{dmath*},
does not converge, the underlying formal result retains its validity at the operatorial level. In other words, the umbral identity remains meaningful, whereas the particular Appell series chosen to realise it is only a local representation and ceases to be admissible at the equilateral point.

Equivalently, one may say that the umbral representation exists globally, while its Appell realisation is restricted to a local convergence chamber.

In a forthcoming investigation, we shall show that a realisation in terms of fractional Franel functionals~\cite{Franel1894,Gould1972}, together with their Mellin--Barnes analytic continuation, provides a natural framework for further developments in this direction.

\section{Declaration of generative AI in the manuscript preparation process}

During the preparation of this work the authors used ChatGPT in order to double-check some of the calculations. After using this tool, the authors reviewed and edited the content as needed and take full responsibility for the content of the published article.

\bibliographystyle{unsrt}
\bibliography{bibliography_bessel_moments}

\end{document}